\documentclass[12pt]{amsart}
\usepackage{amssymb}

\newcommand{\la}{\langle}
\newcommand{\ra}{\rangle}

\newtheorem{theorem}{\bf Theorem}[section]
\newtheorem{lemma}[theorem]{\bf Lemma}

\newtheorem{prop}[theorem]{\bf Proposition}

\newtheorem{example}[theorem]{\bf Example}
\newtheorem{question}[theorem]{\bf Question}

\newcommand{\CC}{{\Bbb C}}

\newcommand{\CP}{{\Bbb CP}}

\newcommand{\QQ}{{\Bbb Q}}
\newcommand{\RR}{{\Bbb R}}

\newcommand{\glie}{{\frak g}}

\newcommand{\klie}{{\frak k}}

\newcommand{\plie}{{\frak p}}

\newcommand{\tlie}{{\frak t}}
\newcommand{\ulie}{{\frak u}}

\newcommand{\solie}{{\frak so}}

\newcommand{\aff}{\operatorname{aff}}
\newcommand{\Ad}{\operatorname{Ad}}
\newcommand{\ad}{\operatorname{ad}}

\newcommand{\End}{\operatorname{End}}

\newcommand{\GL}{\operatorname{GL}}

\newcommand{\Graph}{\operatorname{Graph}}

\newcommand{\Hom}{\operatorname{Hom}}
\newcommand{\Id}{\operatorname{Id}}

\newcommand{\Ker}{\operatorname{Ker}}

\newcommand{\Lie}{\operatorname{Lie}}

\newcommand{\PSL}{\operatorname{PSL}}

\newcommand{\SO}{\operatorname{SO}}

\newcommand{\U}{\operatorname{U}}

\renewcommand{\exp}{\operatorname{exp}}

\newcommand{\ov}{\overline}

\newcommand{\qu}{/\kern-.7ex/}
\newcommand{\exh}{\to\kern-1.8ex\to}

\newcommand{\VP}{{\curly V}\kern-0.9ex\PPP}

\newcommand{\imag}{{\mathbf i}}

\newcommand{\cC}{{\mathcal{C}}}

\newcommand{\fF}{{\mathcal{F}}}
\newcommand{\oO}{{\mathcal{O}}}

\newcommand{\ZII}{Z_{I\kern-.3ex I}}
\newcommand{\ZIII}{Z_{I\kern-.3ex I\kern-.3ex I}}

\date{April 4, 2008}

\author[I. Mundet i Riera]{I. Mundet i Riera}
\address{Departament d'\`Algebra i Geometria, Facultat de Matem\`atiques,
Universitat de Barcelona, Gran Via de les Corts Catalanes 585,
08007 Barcelona, Spain} \email{ignasi.mundet@@ub.edu}

\title{A Hilbert--Mumford criterion for polystability in Kaehler geometry}

\begin{document}
 \maketitle

\begin{abstract}
Consider a Hamiltonian action of a compact Lie group $K$ on a
Kaehler manifold $X$ with moment map $\mu:X\to\klie^*$. Assume
that the action of $K$ preserves the complex structure of $X$, and
consider its unique extension to a holomorphic action of the
complexification $G$ of $K$. We characterize which $G$-orbits in
$X$ intersect $\mu^{-1}(0)$ in terms of the maximal weights
$\lim_{t\to\infty}\la\mu(e^{\imag ts}\cdot x),s\ra$, where
$s\in\klie$. We do not impose any a priori restriction on the
stabilizer of $x$. Under some mild restrictions on the growth of
$\mu$ and the action $K\circlearrowright X$, we view the maximal
weights as defining a collection of maps, for each $x\in X$,
$$\lambda_x:\partial_{\infty}(K\backslash G)\to\RR\cup\{\infty\},$$
where $\partial_{\infty}(K\backslash G)$ is the boundary at
infinity of the symmetric space $K\backslash G$. We prove that
$G\cdot x\cap\mu^{-1}(0)\neq\emptyset$ if: (1) $\lambda_x$ is
everywhere nonnegative, (2) any boundary point $y$ such that
$\lambda_x(y)=0$ can be connected with a geodesic in $K\backslash
G$ to another boundary point $y'$ satisfying $\lambda_x(y')=0$. We
also prove that the maximal weight functions are $G$-equivariant:
for any $g\in G$ and any $y\in
\partial_{\infty}(K\backslash G)$ we have
$\lambda_{g\cdot x}(y)=\lambda_x(y\cdot g)$ .
\end{abstract}

\section{Introduction}

Let $K$ be a compact connected Lie group with Lie algebra $\klie$,
let $X$ be a (non necessarily compact) Kaehler manifold, and
assume that there is an action $K\circlearrowright X$ by
holomorphic isometries (in particular, preserving the symplectic
form) and admitting a moment map $\mu:X\to\klie^*$. Since $\mu$ is
by definition $K$-equivariant, the action of $K$ on $X$ preserves
the zero level set $\mu^{-1}(0)$ and the quotient $\mu^{-1}(0)/K$
carries a natural structure of stratified symplectic manifold, see
\cite{SL}. Let $G$ be the complexification of $K$. By a theorem of
Guillemin and Sternberg \cite{GS}, the action of $K$ on $X$
extends to a unique action of $G$ on $X$ such that the map
$G\times X\ni (g,x)\mapsto g\cdot x\in X$ is holomorphic. The
action of $G$, however, no longer preserves the symplectic form of
$X$ nor the zero level set $\mu^{-1}(0)$.

A very natural and important question, which has been extensively
treated in the literature, is to find a good notion of quotient of
$X$ by the action of $G$, carrying a structure of (possibly
singular) Kaehler manifold induced in some way from the structure
in $X$. The most naive possibility, taking the space of orbits
$X/G$ with the quotient topology, will not even be Hausdorff in
general, so there is no hope to provide it with a structure of
singular Kaehler manifold. To avoid this pathology one can
restrict the attention to a big $G$-invariant subset $X^*\subset
X$, obtained after removing some {\it bad} $G$-orbits in $X$, such
that the quotient space $X^*/G$ is Hausdorff. This faces us with
the problem of defining $X^*$ in a natural way, satisfying the
previous condition and at the same time being as big as possible
(for example, we would like $X^*$ to be dense in $X$).

A systematic way of defining $X^*$ uses the moment map of the
action of $K$. (For the results stated in this paragraph and the
next one, see the survey \cite{HH} and the references therein.)
One says that $x\in X$ is semistable if the closure of $G\cdot x$
intersects $\mu^{-1}(0)$. Then the set $X^{ss}\subset X$ of
semistable points is open and the relation which identifies two
orbits in $X^{ss}$ if their closures intersect is an equivalence
relation. The quotient space $X\qu G$ of this equivalence relation
carries a natural structure of (possibly singular) holomorphic
space, and the projection $X^{ss}\to X\qu G$ is holomorphic. On
the other hand, if one defines $X^{ps}\subset X$ as the set of
points whose $G$-orbit intersects $\mu^{-1}(0)$ (we call such
points polystable), then the inclusion $X^{ps}\subset X^{ss}$
induces a homeomorphism $X^{ps}/G\simeq X\qu G$, so that one can
take $X^{ps}$ as a good choice for $X^*$. This motivates the
following question.

\begin{question}
\label{question:2} Which $G$-orbits $\oO=G\cdot x\subset X$
intersect $\mu^{-1}(0)$? If $\oO$ is such an orbit, how many
$K$-orbits does $\oO\cap\mu^{-1}(0)$ contain?
\end{question}

In the first question we would like some characterization of the
points $x\in X$ such that $G\cdot x\cap \mu^{-1}(0)\neq\emptyset$
in terms of the symplectic geometry of the action
$K\circlearrowright X$. The answer to the second question (namely,
that $\oO\cap\mu^{-1}(0)$ contains at most one $K$-orbit) is well
known and follows from an easy argument which we recall in Section
\ref{ss:proof-1-2} below. A refinement of this is the statement
that $X\qu G$ is homeomorphic to the symplectic quotient
$\mu^{-1}(0)/K$, and that the holomorphic structure on $X\qu G$
defined in \cite{HH} is compatible with it, in the sense that they
combine to define a structure of stratified Kaehler manifold (see
\cite{HH,S}).

A partial characterization of the $G$-orbits intersecting
$\mu^{-1}(0)$ was given in \cite{M} in terms of the maximal
weights $\lambda(x;s)$, defined for any $s\in\klie$ to be
$$\lambda(x;s)=\lim_{t\to\infty}\la\mu(e^{\imag ts}\cdot
x),s\ra\in\RR\cup\{\infty\}$$ (this limit exists, see Section
\ref{ss:asymptotics} below). A point $x\in X$ was defined to be
analytically stable if $\lambda(x;s)>0$ for any nonzero
$s\in\klie$, and  \cite[Theorem 5.4]{M} states that $x$ is
analytically stable if and only if $G\cdot
x\cap\mu^{-1}(0)\neq\emptyset$ and the stabilizer $G_x$ is finite.
Apart from the restriction to points with finite stabilizers, the
point of view in \cite{M} has the inconvenient that it does not
seem to allow any straightforward proof that if $x$ is
analytically stable then $g\cdot x$ is also analytically stable
for any $g\in G$. Roughly speaking, this is proved in \cite{M} as
a consequence of the characterization of analytic stability in
terms of the so-called linear properness of the integral of the
moment map, which essentially amounts to \cite[Theorem 5.4]{M}.
(The equivariance property of the moment map allows to prove that
$\lambda(k\cdot x;s)=\lambda(x;\Ad(k)(s))$ for any $k\in K$, which
clearly implies that analytic stability is a property of
$K$-orbits, but there is no obvious action of $G$ on $\klie$
extending the adjoint action of $K$ giving a similar
$G$-equivariance property of the maximal weights.)

It is easy to deduce from the results in \cite{M} that if
$\lambda(x;s)<0$ for some $s$ then $G\cdot
x\cap\mu^{-1}(0)=\emptyset$. Hence if $G\cdot x\cap\mu^{-1}(0)\neq
\emptyset$ then $\lambda(x;s)\geq 0$ for any $s$. But this is not
a sufficient condition: to decide whether $G\cdot
x\cap\mu^{-1}(0)$ is nonempty one needs to control in some way
which of the maximal weights vanish.

A well known and elementary example which illustrates these
phenomena is given by the action of the isometries of $S^2$ on
tuples of points. Here $S^2$ denotes the unit sphere in $\RR^3$
centered at $0$ and endowed with the round metric. Let
$K=\SO(3,\RR)$ be the group of orientation preserving isometries
of $S^2$. The complexification of $K$ is $G=\PSL(2,\CC)$, which
can be identified with the holomorphic automorphisms of
$S^2\simeq\CP^1$. Let $X=(S^2)^4$ and take on $X$ the product
Kaehler structure. The diagonal action of $K$ on $X$ clearly
preserves the holomorphic structure and the symplectic form, and a
moment map for it is given by sending any tuple
$(x_1,\dots,x_4)\in S^2$ to its center of mass
$\frac{1}{4}(x_1+\dots+x_4)$ (here we identify
$\RR^3\simeq\solie(3,\RR)^*$ using the vector product in $\RR^3$).
The reader can easily check that, if $\{x_1,x_2,x_3,x_4\}\subset
S^2$ are distinct points, then
\begin{itemize}
\item $x=(x_1,x_2,x_3,x_4)\in X$ is analytically semistable, \item
$x'=(x_1,x_1,x_2,x_3)\in X$ satisfies $G\cdot
x'\cap\mu^{-1}(0)=\emptyset$ but all maximal weights
$\lambda(x';s)$ are nonnegative, \item $x''=(x_1,x_1,x_2,x_2)\in
X$ satisfies $G\cdot x''\cap\mu^{-1}(0)\neq\emptyset$, but some of
the maximal weights $\lambda(x;s)$ vanish.
\end{itemize}
(Of course, $x'$ is semistable in the usual sense in GIT, and the
closure of orbit $G\cdot x'$ contains $G\cdot x''$ and hence meets
$\mu^{-1}(0)$.)

A complete solution to Question \ref{question:2} was given by A.
Teleman in \cite{T}. However, the result in \cite{T} has some
limitations. First, it is assumed that $X$ satisfies a condition
called energy-completeness (see \cite[Definition 3.8]{T}). Second,
when giving a sufficient condition for a point $x\in X$ to satisfy
$G\cdot x\cap\mu^{-1}(0)$ it is assumed that the Lie algebra
$\glie_x$ of the stabilizer $G_x$ is reductive (see
\cite[Definition 3.12]{T}). This is a little bit unsatisfactory:
it might be preferable to obtain the reductivity of $\glie_x$ as a
consequence of a simpler condition involving exclusively the
maximal weights, not any information on the stabilizer of $x$.

In this paper we propose an alternative answer to Question
\ref{question:2} based on viewing the maximal weights as defining
a function on the boundary $\partial_{\infty}(K\backslash G)$ of
the symmetric space $K\backslash G$. Such boundary exists by the
general theory of Hadamard spaces, of which $K\backslash G$ is an
instance (see \cite{B,E}), and it is homeomorphic to a sphere of
dimension one unit less than that of $K\backslash G$. To prove our
results we still require some technical restrictions to be
satisfied by $X$; namely, we assume that the moment map (resp. the
vector fields generated by the infinitesimal action) grows
quadratically (resp. linearly) with respect to the distance
function from a given base point. More precisely: given a
biinvariant metric on $\klie$ we require that there exists a point
$x_0\in X$ and a constant $C$ such that for any $x\in X$ and any
$s\in\klie$ we have
\begin{align}
|\xi_s(x)| & \leq C\,|s|\,(1+d_X(x,x_0)) \label{eq:bound-action} \\
|\mu(x)| & \leq C\,(1+d_X(x,x_0)^2) \label{eq:bound-moment-map}
\end{align}
where $\xi_s\in\cC^{\infty}(TX)$ is the vector field generated by
the infinitesimal action of $s$ and $d_X$ is the distance function
between points in $X$. These conditions are satisfied e.g. when
$X$ is compact or when $X$ is a vector space and the action of $K$
is linear. (On the other hand, in this paper we do not assume any
completeness condition as in \cite{T}.)

Assuming (\ref{eq:bound-action}) and (\ref{eq:bound-moment-map})
we construct in Section \ref{ss:extending-integral} the maximal
weight function
$$\lambda_x:\partial_{\infty}(K\backslash G)\to\RR\cup\{\infty\}$$
for any $x\in X$ as an appropriate limit of a normalization of the
integral of the moment map
$$\psi_x:K\backslash G\to\RR.$$
The integral of the moment map was defined in \cite{M} as a
function $G\to\RR$, and it was observed in \cite[Proposition
3.4]{M} that it is invariant under the action of $K$ on the left
on $G$. The boundary $\partial_{\infty}(K\backslash G)$ carries an
action of $G$ extending the right action on $K\backslash G$ by
isometries, and we prove in Lemma \ref{lemma:lambda-equivariant}
that for any $x\in X$ and $g\in G$ we have
\begin{equation}
\label{eq:lambda-equivariant} \lambda_{g\cdot
x}(y)=\lambda_x(y\cdot g). \end{equation} This property is a
consequence of the cocycle property satisfied by the integral of
the moment map (see formula (\ref{eq:cocycle}) below).

We say that $x\in X$ is analytically stable if for any $y\in
\partial_{\infty}(K\backslash G)$ we
have $\lambda_x(y)>0$. A point $x\in X$ is said to be analytically
polystable if for any $y\in
\partial_{\infty}(K\backslash G)$ we have $\lambda_x(y)\geq 0$
and for any $y\in\partial_{\infty}(K\backslash G)$ such that
$\lambda_x(y)=0$ there exists some $y'\in
\partial_{\infty}(K\backslash G)$
such that $\lambda_x(y')=0$ and the points $y,y'$ can be connected
by a geodesic in $K\backslash G$.

The next theorem is the main result of the paper. It will be
proved in Section \ref{s:proof-theorem}.

\begin{theorem}
\label{thm:analytic-stab-zero-moment-map} Let $x\in X$ be any
point and let $G_x=\{g\in G\mid g\cdot x=x\}$ be its stabilizer.
\begin{enumerate}
\item If $x$ is analytically stable (resp. analytically
polystable) then $g\cdot x$ is analytically stable (resp.
analytically polystable) for each $g\in G$. \item The intersection
$G\cdot x\cap\mu^{-1}(0)$ consists of at most one $K$-orbit. \item
$x$ is analytically stable if and only if $G\cdot
x\cap\mu^{-1}(0)\neq\emptyset$ and $G_x$ is finite. \item $x$ is
analytically polystable if and only if $G\cdot
x\cap\mu^{-1}(0)\neq\emptyset$. If this is the case, then $G_x$ is
reductive.
\end{enumerate}
\end{theorem}

Note that statement (1) follows immediately from the equivariance
property (\ref{eq:lambda-equivariant}). (2) and (3) are well
known, but we also prove them for completeness (the proof we give
of (3), using the index of maps between spheres, is new to the
best of our knowledge).

For any $s\in\klie$ of unit length, let $e_s$ be the boundary
point of $K\backslash G$ to which the geodesic $t\mapsto[e^{\imag
ts}]$ converges as $t\to\infty$. We have:
$\lambda(x;s)=\lambda_x(e_s)$ for any $s\in\klie$ of unit norm,
$\lambda(x;\alpha s)=\alpha\lambda(x;s)$ for any positive real
number $\alpha$, and any point in $\partial_{\infty}(K\backslash
G)$ is of the form $e_s$ for some $s$. Hence, the notion of
analytic stability in the present paper coincides with the notion
given in \cite{M}. In order to rephrase the definition of
polystability in terms of the functions $\lambda(x;s)$ we
introduce the following definitions. Two skew Hermitian
endomorphisms $a,b$ of a complex vector space $V$ are said to be
opposed if $\imag a$ and $-\imag b$ have the same spectrum, say
$\{\lambda_0<\dots<\lambda_r\}\subset\RR$, and the growing
filtrations $W_a^{\bullet}$ and $W_b^{\bullet}$ defined as
$$W_a^j=\bigoplus_{i\leq j}\Ker(\imag a-\lambda_i\Id),
\qquad\qquad W_b^j=\bigoplus_{i\geq r-j}\Ker(\imag
b+\lambda_i\Id)$$ satisfy
$$V=\bigoplus_{p+q=r}W_a^p\cap W_b^q.$$
Two elements $u,v\in\klie$ are said to be opposed if $u$ and $-v$
belong to the same adjoint orbit in $\klie$ and $\ad(u),\ad(v)$
are opposed endomorphisms of $\glie$. For example, for any
$u\in\klie$ the elements $u,-u$ are opposed. Note also that if $u$
has unit norm and $-v$ belongs to the adjoint orbit through $u$,
then $v$ also has unit norm because the norm in $\klie$ is
biinvariant. The following lemma is a consequence of Lemma
\ref{lemma:opposed-geodesically-connected} below.

\begin{lemma}
\label{lemma:rephrasing-polystability} A point $x\in X$ is
polystable if $\lambda(x;s)\geq 0$ for any $s\in\klie$ and if, for
any nonzero $s\in\klie$ such that $\lambda(x;s)=0$, there exists
some $u\in\klie$ which is opposed to $s$ and such that
$\lambda(x;u)=0$.
\end{lemma}

The results in this paper can be seen as an analytic version of
part of the results in Chapter 2 of \cite{MFK}. Mumford's point of
view is that the maximal weights in the case of projective
varieties define a function on the flag complex $\Delta(G)$, which
is the set of rational points at infinity of $G$ and can be
naturally thought as an algebraic version of the boundary
$\partial_{\infty}(K\backslash G)$. More precisely, the function
$\Delta(G)\ni\delta\mapsto\nu^L(x,\delta)$ defined in p. 59 of
[op. cit.] is the analogue of our function $\lambda_x$. When $X$
is projective and its symplectic structure is the restriction of
the Fubini--Study structure on the projective space, statement (3)
in Theorem \ref{thm:analytic-stab-zero-moment-map}, combined with
Kempf--Ness's results (see \cite{Sch} for an excellent survey),
implies the usual Hilbert--Mumford criterion for stability in GIT,
and this explains the title of the present paper. On the other
hand, statement (4) in Theorem
\ref{thm:analytic-stab-zero-moment-map} gives a characterization
of which points $x$ in a linear representation of a reductive
group $G$ have closed orbit $G\cdot x$ in terms uniquely of
maximal weights, and this seems to be a new result (note that
antipodal points in $\Delta(G)$, as defined in Definition 2.8, p.
61 in [op. cit.], correspond to pairs of points in
$\partial_{\infty}(K\backslash G)$ which can be connected by a
geodesic).

The rest of the paper is organized as follows. In Section
\ref{s:boundary-K-G} we recall the definition and some basic facts
on the boundary at infinity of the symmetric space $K\backslash
G$. In Section \ref{s:maximal-weights} we construct the maximal
weight functions $\lambda_x:\partial_{\infty}(K\backslash
G)\to\RR\cup\{\infty\}$. In Section \ref{s:proof-theorem} we give
the proof of Theorem \ref{thm:analytic-stab-zero-moment-map} and,
finally, in Section
\ref{s:proof:lemma:opposed-geodesically-connected} we prove Lemma
\ref{lemma:opposed-geodesically-connected}.

\section{The symmetric space $K\backslash G$ and its boundary
at infinity $\partial_{\infty}(K\backslash G)$}
\label{s:boundary-K-G}

\subsection{The boundary $\partial_{\infty}(K\backslash G)$}
The coset space $K\backslash G$ has a natural structure of
differentiable manifold. We consider on it the action of $G$ given
by multiplication on the right: $[g]\cdot h=[gh]$ for any $g,h\in
G$. Let $x_0\in K\backslash G$ denote the class of the identity
element $1_G\in G$. Choose a biinvariant Euclidean norm on
$\klie$. This induces a unique $G$-invariant Riemannian metric on
$K\backslash G$, because the action $T(K\backslash
G)\circlearrowleft G$ given by differentiating right
multiplication is transitive, the stabilizer of the fiber
$T_{x_0}(K\backslash G)$ over the identity element $1_G\in G$ is
$K$, and the action of $K$ on $T_{x_0}(K\backslash G)$ can be
identified with the adjoint action of $K$ on $\klie$ (via the
natural identification $T_{x_0}(K\backslash G)\simeq\imag\klie$).
The geodesics corresponding to this metric are given by maps
$t\mapsto [e^{\imag t s}g]\in K\backslash G$ for any $s\in\klie$
and $g\in G$.

The invariant metric on $K\backslash G$ has nonpositive curvature
(see \cite{E}) so, endowed with it, $K\backslash G$ is a Hadamard
space. So the general theory of Hadamard spaces (see for example
\cite{B}) implies that there is a naturally defined boundary at
infinity $\partial_{\infty}(K\backslash G)$. This can be described
in concrete terms using geodesic rays i.e. maps
$$\gamma:(0,\infty)\to K\backslash G$$ giving a parametrization by
arc of a portion of geodesic. Let $d$ denote the distance function
between points in $K\backslash G$. Two geodesic rays
$\gamma_0,\gamma_1$ are declared to be equivalent
$\gamma_0\sim\gamma_1$ if the distance
$d(\gamma_0(t),\gamma_1(t))$ is bounded independently of $t$. This
is an equivalence relation on the set of geodesic rays, and the
boundary at infinity of $K\backslash G$ is the set of equivalence
classes:
$$\partial_{\infty}(K\backslash G)=\{\text{ geodesic rays
}\}/\sim.$$ If $\gamma:(0,\infty)\to K\backslash G$ is a geodesic
ray and $g\in G$ then we define $\gamma\cdot g$ to be the geodesic
ray whose value at $t$ is $\gamma(t)\cdot g$. This defines a right
action of $G$ on the set of geodesic rays. Since the action of $G$
on the right on $K\backslash G$ is by isometries, this action on
the set of geodesic rays preserves the equivalence $\sim$ and
hence descends to an action on $\partial_{\infty}(K\backslash G)$.

Let $S(\klie)\subset\klie$ denote the unit sphere. For any $s\in
S(\klie)$ we define $e_s\in\partial_{\infty}(K\backslash G)$ to be
the class of the geodesic ray $\eta_s:(0,\infty)\to K\backslash G$
defined as $\eta_s(t)=[e^{\imag ts}]$. Then the map $e:S(\klie)\ni
s\mapsto [e_s]\in\partial_{\infty}(K\backslash G)$ is a bijection
(see Section II.2 in \cite{B}). We endow
$\partial_{\infty}(K\backslash G)$ with the topology which makes
$e$ a homeomorphism. Then the action of $G$ on
$\partial_{\infty}(K\backslash G)$ is by homeomorphisms. For each
$s\in S(\klie)$ and any $g\in G$ define $s\cdot g\in S(\klie)$ by
the property that
$$e_s\cdot g=e_{s\cdot g}.$$

We remark that the boundary $\partial_{\infty}(K\backslash G)$ is
independent of the chosen biinvariant metric on $\klie$. Indeed,
geodesic rays do not depend on the choice of metric (they are
always of the form $t\mapsto [e^{\imag s t}g]$) and neither does
the equivalence relation $\sim$ on geodesic rays, because the
distance functions on $K\backslash G$ induced by two choices of
biinvariant metric on $\klie$ are uniformly comparable.

\subsection{The case $K=\U(n)$ and $G=\GL(n,\CC)$}
\label{ss:boundary-unitary} When $K=\U(n)$ (so that
$G=\GL(n,\CC)$) the action of $G$ on $S(\ulie(n))$ can be computed
in concrete terms, as we will shortly see. Define the logarithm
map $\log:G\to\klie$ by the condition that $\log(g)=u$ if $g=k
e^{\imag u}$ is the Cartan decomposition of $g$, so that $k\in K$
and $u\in\klie$. Let $s\in S(\ulie(n))$. The matrix $\imag s$ is
Hermitian symmetric, so it diagonalizes and has real eigenvalues,
say $\lambda_1<\dots<\lambda_r$. Let $V_j=\Ker(\lambda_j-\imag s)$
be the eigenspace corresponding to $\lambda_j$ and define
$V^k=V_1\oplus\dots\oplus V_k$ for any integer $k\geq 1$. Take any
$g\in G$ and define
$$V_j^{\infty}=(g^{-1}(V_{j-1}))^{\perp}\cap g^{-1}(V_j),$$
where $V^{\perp}$ denotes the orthogonal of $V$. Then we have a
direct sum decomposition $\CC^n=\bigoplus V_j^{\infty}$. Define
$\rho_g(s)\in\ulie(n)$ by the conditions that $\rho_g(s)$
preserves each $V_j^{\infty}$ and that the restriction of
$\rho_g(s)$ to $V_j^{\infty}$ is given by multiplication by
$-\imag\lambda_j$. We claim that $\rho_g(s)$ is equal to $s\cdot
g$. This is equivalent to the statement
\begin{equation}
\label{eq:computing-rho-g}
\rho_g(s)=\lim_{\tau\to\infty}\frac{1}{\tau}\log(e^{\imag \tau
s}g).
\end{equation} To prove (\ref{eq:computing-rho-g}) one can argue
as follows. Take a very small $\epsilon>0$ (in particular, smaller
than $\inf\{\lambda_j-\lambda_{j-1}\}/3$). Using the variational
description of eigenvalues and eigenspaces of $\log(h)$, one
proves that for big enough $\tau$ the eigenvalues of
$s_{\tau}:=\tau^{-1}\log(e^{\imag \tau s}g)$ are contained in
$\bigcup [\lambda_j-\epsilon,\lambda_j+\epsilon]$, and the number
of eigenvalues in $[\lambda_j-\epsilon,\lambda_j+\epsilon]$ is
equal to $\dim V_j$. Let $V_j^{\tau}$ be the direct sum of the
eigenspaces of $s_{\tau}$ with eigenvalue contained in
$[\lambda_j-\epsilon,\lambda_j+\epsilon]$. Then $V_j^{\tau}$
converges to $V_j^{\infty}$ in the Grassmannian variety. Details
are left to the reader.

Using the previous computations, we can also check that the map
$e:S(\klie)\to\partial_{\infty}(K\backslash G)$ is a bijection.
This is equivalent to proving the existence of a bound, for each
$g\in G$ and $s\in S(\klie)$, of the form $d([e^{\imag t (s\cdot
g)}],[e^{\imag t s}g])\leq C,$ where the constant $C>0$ is
independent of $t$. Indeed, this implies that the geodesic ray
$t\mapsto [e^{\imag t s}g]$ is equivalent to $t\mapsto[e^{\imag
t(s\cdot g)}]$, which is $e_{s\cdot g}$. Details are left as an
exercise to the reader.

\subsection{Tori generated by elements in $\klie$}
For any $s\in\klie$ we define the torus
$$T_s=\ov{\{\exp(ts)\mid t\in\RR\}}\subset K.$$

\begin{lemma}
\label{lemma:dim-torus-invariant} For any $s\in\klie$ and any
$g\in G$ we have $\dim T_s=\dim T_{s\cdot g}$.
\end{lemma}
\proof
By Peter--Weyl theorem one can pick an embedding of Lie
groups $K\hookrightarrow\U(n)$ which complexifying induces an
inclusion $G\hookrightarrow\GL(n,\CC)$. Since the boundary at
infinity does not depend on the choice of biinvariant metric, this
inclusion induces an inclusion of boundaries
$\partial_{\infty}(K\backslash G)\hookrightarrow
\partial_{\infty}(U(n)\backslash\GL(n,\CC))$, which is equivariant
with respect to the natural action of $G$ on
$\partial_{\infty}(K\backslash G)$ and the action of $G$ on
$\partial_{\infty}(U(n)\backslash\GL(n,\CC))$ given by the
inclusion $G\hookrightarrow\GL(n,\CC)$ (see the proof of Lemma
\ref{lemma:phi-G-equivariant} for details). All this implies that
it suffices to consider the case $K=\U(n)$. But if $s\in\ulie(n)$
then the dimension of $T_s$ depends uniquely on the eigenvalues of
$s$ (namely, $\dim T_s$ is equal to the dimension of the
$\QQ$-vector space spanned by the eigenvalues of $s$). On the
other hand, the observations in Section \ref{ss:boundary-unitary}
imply that for any $s\in\ulie(n)$ and $g\in \GL(n,\CC)$ the
endomorphisms $s,s\cdot g\in\End\CC^n$ have the same set of
eigenvalues, so we certainly have $\dim T_s=\dim T_{s\cdot g}$.
\qed

\subsection{Geodesically connected points}
Two points in $\partial_{\infty}(K\backslash G)$ are said to be
geodesically connected if there is a geodesic in $K\backslash G$
which converges on one side to one of the points and on the other
side to the other point. This definition is independent of the
biinvariant metric on $\klie$ because the set of geodesics in
$K\backslash G$ and the notion of convergence of rays to points in
$\partial_{\infty}(K\backslash G)$ do not depend on the metric on
$\klie$. A trivial example:

\begin{example}
\label{example:s-minus-s} For any $s\in S(\klie)$ the points $e_s,
e_{-s}\in\partial_{\infty}(K\backslash G)$ are geodesically
connected by the geodesic $t\mapsto [e^{\imag ts}]$.
\end{example}

A concrete translation into algebraic terms of the condition of
being geodesically connected can be given using the notion of
opposed elements in $\klie$ defined in the Introduction:

\begin{lemma}
\label{lemma:opposed-geodesically-connected} Given $u,v\in
S(\klie)$, the points $e_u,e_v\in\partial_{\infty}(K\backslash G)$
are geodesically connected if and only if $u,v$ are opposed.
\end{lemma}

To avoid an excessive detour from our arguments we postpone the
proof of the lemma to Section
\ref{s:proof:lemma:opposed-geodesically-connected} at the end of
the paper.

A more synthetic characterization of geodesic connectedness may be
given in terms of parabolic subgroups. We state such translation
for the sake of completeness, but we will not use it in the
sequel. Recall that a parabolic subgroup of $G$ is by definition
the stabilizer of a point in $\partial_{\infty}(K\backslash G)$.
It is almost a tautology that the stabilizer of
$e_s\in\partial_{\infty}(K\backslash G)$ is the subgroup
$P_s\subset G$ consisting of all $g\in G$ such that $e^{\imag t
s}g e^{-\imag t s}$ stays bounded as $t\to\infty$, so that all
parabolic subgroups of $G$ are of the form $P_s$ for some $s\in
S(\klie)$. The maximal reductive subgroups of $P_s$ are called the
Levi subgroups (they are all pairwise conjugate). Two parabolic
subgroups $P_s,P_{s'}\subset G$ are said to be opposed if $P_s\cap
P_{s'}$ is a Levi subgroup both of $P_s$ and $P_{s'}$. Now, $e_s$
and $e_{s'}$ are geodesically connected if and only if $P_s$ and
$P_{s'}$ are opposed and $s,-s'$ belong to the same coadjoint
orbit in $\klie$.

\section{Maximal weights as a map
$\lambda_x:\partial_{\infty}(K\backslash G)\to \RR\cup\{\infty\}$}
\label{s:maximal-weights}

We now come back to the situation considered in the Introduction,
so that $K\circlearrowright X$ is a Hamiltonian action of a
compact Lie group $K$ on a Kaehler manifold $X$ preserving the
complex structure, and we consider the extension of this action to
a holomorphic action $G\circlearrowright X$ of the
complexification $G=K^{\CC}$.

\subsection{The integral of the moment map} Denote by
$\pi:\glie=\klie\oplus\imag\klie\to\imag\klie$ the projection to
the second factor. Let $r_{g^{-1}}:G\to G$ be the map given by
multiplication by $g^{-1}$ on the right, and let
$Dr_{g^{-1}}:T_gG\to T_{1_G}G\simeq\glie$ be its derivative. For
any $v\in T_gG$ we define $v\cdot g^{-1}:=Dr_{g^{-1}}(v)\in\glie$.

For any $x\in X$ we define a one form $\sigma_x\in\Omega^1(G)$ as
follows:
$$\sigma_x(g)(v):=\la\mu(g\cdot x),-\imag\pi(v\cdot g^{-1})\ra$$
for any $v\in T_vG$. It is immediate to deduce from the definition
that for any $g,h\in G$ and any $v\in T_gG$ we have
$\sigma_x(gh)(v\cdot h)=\sigma_x(g)(v)$, so that
\begin{equation}
\label{eq:sigma-equivariant} \sigma_{hx}=r_h^*\sigma_x.
\end{equation}
By \cite[Lemma 3.1]{M} the form $\sigma_x$ is exact. Hence we may
define $\Psi_x:G\to\RR$ to be the unique function such that
$\Psi_x(1_G)=0$ and $d\Psi_x=\sigma_x$. We call $\Psi_x$ the {\bf
integral of the moment map}. Property (\ref{eq:sigma-equivariant})
implies the following cocycle formula:
\begin{equation}
\label{eq:cocycle} \Psi_x(g)+\Psi_{g\cdot x}(h)=\Psi_x(hg)
\end{equation}
for any $x\in X$ and $g,h\in G$.

\subsection{Asymptotics of the integral of the moment map}
\label{ss:asymptotics} Given $s\in\klie$ we define
$\mu_s(x)=\la\mu(x),s\ra$ for any $x\in X$ and for any $t\in\RR$
we define $\lambda_t(x,s)=\mu_s(e^{\imag t s}\cdot x)$. For any
$s\in\glie$ we denote by
$$\xi_s\in\cC^{\infty}(TX)$$
the vector field generated by the infinitesimal action of $s$.
Since the action of $G$ on $X$ is holomorphic we have $\xi_{\imag
s}=I\xi_s$, where $I\in\cC^{\infty}(\End TX)$ is the complex
structure on $X$. Using the defining properties of the moment map
we compute:
\begin{align}
\partial_t\lambda_t(x;s) &=
\partial_t\la\mu(e^{\imag t s}\cdot x),s\ra
=\omega(\xi_s,I\xi_s)(e^{\imag t s}\cdot x) \notag \\
\label{eq:gradient-mu-s} &=\la \xi_s(e^{\imag t s}\cdot
x),\xi_{u}(e^{\imag t s}\cdot x)\ra= |I\xi_s(e^{\imag t s}\cdot
x)|^2,
\end{align} where
$\partial_t$ denotes the derivative with respect to $t$. This
implies
\begin{equation}
\label{eq:integral-t-pes} \lambda_t(x;s)=\la\mu(x),s\ra+\int_0^t
|\xi_s(e^{\imag \tau s}\cdot x)|^2\,d\tau,
\end{equation}
and in particular $\lambda_t(x;s)$ is nondecreasing as a function
of $t$.

It follows from the definition of $\sigma_x$ that for any
$s\in\klie$
$$\Psi_x(e^{\imag ts})=\int_0^{t}\lambda_{\tau}(x;s)\,d\tau.$$
Since $\lambda_{\tau}$ is nondecreasing we deduce that
\begin{equation}
\label{eq:limit-integral-mom-map}
\lim_{t\to\infty}\frac{\Psi_x(e^{\imag ts})}{t}=\lambda(x;s):=
\lim_{t\to\infty}\lambda_t(x;s)\in\RR\cup\{\infty\}.
\end{equation}
The limit $\lambda(x;s)$ is what was defined to be the maximal
weight in \cite{M}. When it is necessary to be more specific, we
will say that $\lambda(x;s)$ is the maximal weight of the action
of $K$ on $X$ and we will denote it by $\lambda^K(x;s)$ (this will
be the case in Section \ref{s:proof-theorem}, where different
symmetry groups will be considered simultaneously).

We end this section by showing how the growth of the integral of
the moment map can be used to bound the distance between points in
$X$. Recall that $d_X$ denotes the distance function between pairs
of points in $X$.

\begin{lemma}
\label{lemma:integral-controls-distance} Let $g\in G$ and
$s\in\klie$. If $\Psi_x(e^{\imag ts}g)t^{-1}$ is bounded uniformly
on $t$, then $d_X(e^{\imag ts}g\cdot x,x)t^{-1/2}$ converges to
$0$ as $t\to\infty$.
\end{lemma}
\proof  Using (\ref{eq:cocycle}) and
(\ref{eq:limit-integral-mom-map}) we compute
$$\lim_{t\to\infty}\frac{\Psi_x(e^{\imag ts}g)}{t}=
\lim_{t\to\infty}\frac{\Psi_x(g)+\Psi_{g\cdot x}(e^{\imag ts})}{t}
=\lim_{t\to\infty}\frac{\Psi_{g\cdot x}(e^{\imag
ts})}{t}=\lambda(g\cdot x;s).$$ So, if $\Psi_x(e^{\imag
ts}g)t^{-1}$ is bounded uniformly on $t$, then, by
(\ref{eq:integral-t-pes}), $\int_0^{\infty} |\xi_s(e^{\imag \tau
s}g\cdot x)|^2\,d\tau<\infty.$ Since, on the other hand,
$$d_X(e^{\imag t s}g\cdot x,g\cdot x)\leq\int_0^t
|\xi_{\imag s}(e^{\imag \tau s}g\cdot x)|\,d\tau =\int_0^t
|\xi_{s}(e^{\imag \tau s}g\cdot x)|\,d\tau,$$ the following lemma
applied to $f(\tau)=|\xi_s(e^{\imag \tau s}g\cdot x)|$ implies
that $d_X(e^{\imag ts}g\cdot x,g\cdot x)t^{-1/2}$ converges to $0$
as $t\to\infty$. The lemma is finished by applying the triangular
inequality. \qed

\begin{lemma}
Let $f:(0,\infty)\to\RR_{\geq 0}$ be a nonnegative square
integrable function, so that we have
$\int_0^{\infty}f^2(\tau)\,d\tau<\infty$. Then
$$\left(\int_0^tf(\tau)\,d\tau\right)t^{-1/2}\to 0\qquad\text{as
$t\to\infty$.}$$
\end{lemma}
\proof Let $E=\int_0^{\infty}f^2(\tau)\,d\tau$, let $\epsilon>0$
be any real number and choose $t_0>0$ in such a way that
$\int_0^{t_0}f(\tau)^2\,d\tau\geq (1-\epsilon)E$, so that for any
$t\geq t_0$ we have $$\int_{t_0}^tf(\tau)^2\,d\tau\leq \epsilon
E.$$ Then we compute, using Cauchy--Schwartz and the previous
estimate:
\begin{align*}
\int_0^tf(\tau)\,d\tau &=
\int_0^{t_0}f(\tau)\,d\tau+\int_{t_0}^tf(\tau)\,d\tau =
\int_0^{t_0}f(\tau)\,d\tau+\int_{t_0}^t(\epsilon^{-1/2}f(\tau))
\epsilon^{1/2}\,d\tau \\
&\leq (E t_0)^{1/2}+(E \epsilon (t-t_0))^{1/2}
=E^{1/2}(t_0^{1/2}+\epsilon^{1/2}(t-t_0)^{1/2}).
\end{align*}
The result follows by observing that
$(t_0^{1/2}+\epsilon^{1/2}(t-t_0)^{1/2})t^{-1/2}\to\epsilon^{1/2}$
as $t\to\infty$. \qed

\subsection{Extending the integral of the moment map
to $\partial_{\infty}(K\backslash G)$}
\label{ss:extending-integral} The cocycle condition
(\ref{eq:cocycle}) and the fact that for any $y\in X$ the
restriction of $\Psi_y$ to $K\subset G$ vanishes identically
(which follows immediately from the definition) implies that
$\Psi_x(kg)=\Psi_x(g)$ for each $k\in K$ and $g\in G$, so that
$\Psi_x$ descends to a map
$$\psi_x:K\backslash G\to\RR.$$ Recall that $d$ denotes the distance
function between pairs of points in $K\backslash G$, and that
$x_0\in K\backslash G$ denotes the class of the identity element
in $G$. We are next going to prove that the function
$\phi_x:K\backslash G\to \RR$ defined as
$$\phi_x(z)=\frac{\psi_x(z)}{d(z,x_0)}$$ extends to a function
on the boundary $\partial_{\infty}(K\backslash G)$. For any
geodesic ray $\gamma:[0,\infty)\to K\backslash G$ we define
$$\lambda_x(\gamma)=
\lim_{t\to\infty}\phi_x(\gamma(t))\in\RR\cup\{\infty\}.$$ The
extendability of $\phi_x$ is equivalent to the following lemma.

\begin{prop} If the geodesic rays $\gamma_0,\gamma_1$ satisfy
$\gamma_0\sim\gamma_1$, then
$$\lambda_x(\gamma_0)=\lambda_x(\gamma_1).$$
\end{prop}
\proof We may write $\gamma_j(t)=[e^{\imag ts_j}g_j]$ for $j=0,1$,
where $s_j\in\klie$ and $g_j\in G$. Assume that
$\lambda_x(\gamma_0)$ is finite. Then Lemma
\ref{lemma:integral-controls-distance} implies that $d(e^{\imag
ts_0}g_0\cdot x,x_0)t^{-1/2}$ converges to $0$ as $t\to\infty$.
Since $\gamma_0\sim\gamma_1$, we may bound
$d(\gamma_0(t),\gamma_1(t))\leq \kappa$ uniformly on $t$. It
follows that, for any $t$, we may take a smooth function
$\rho_t:[0,1]\to G$ such that $\rho_t(j)=e^{\imag s_jt}g_j$ for
$j=0,1$, and such that
\begin{equation}
\label{eq:cami-curt} \int_0^1|\partial_{\nu}\rho_t(\nu)\cdot
\rho_t(\nu)^{-1}|\,d\nu\leq \kappa.
\end{equation}
This bound, together with $d(e^{\imag ts_0}g_0\cdot
x,x_0)t^{-1/2}\to 0$ and assumption (\ref{eq:bound-action}),
implies that, for any $\nu$, $d(\rho_t(\nu)\cdot x,x_0)t^{-1/2}\to
0.$ Using assumption (\ref{eq:bound-moment-map}) we get
$|\mu(\rho_t(\nu)\cdot x)|t^{-1}\to 0.$ Since $d\Psi_x=\sigma_x$,
the previous formula together with (\ref{eq:cami-curt}) implies
that
$$\lim_{t\to\infty}\frac{|\Psi_x(e^{\imag ts_0}g_0)-\Psi_x(e^{\imag
ts_1}g_1)|}{t}=0,$$ from which we deduce
$\lambda_x(\gamma_0)=\lambda_x(\gamma_1)$. This immediately
implies, arguing by contradiction, that if
$\lambda_x(\gamma_0)=\infty$ then $\lambda_x(\gamma_1)=\infty$.
For if $\lambda_x(\gamma_1)<\infty$ then, reversing the roles of
$\gamma_0$ and $\gamma_1$ in the previous arguments, we deduce
that $\lambda_x(\gamma_0)=\lambda_x(\gamma_1)<\infty$. \qed

Using the previous lemma we may define
$\lambda_x(y):=\lambda_x(\gamma)\in\RR\cup\{\infty\}$ for any
$y\in\partial_{\infty}(K\backslash G)$, where $\gamma$ is any
geodesic ray representing $y$. In this way we obtain a well
defined map
$$\lambda_x:\partial_{\infty}(K\backslash
G)\to\RR\cup\{\infty\},$$ which we call the maximal weight
function. We now prove a crucial equivariance property of the
maximal weights.

\begin{lemma}
\label{lemma:lambda-equivariant} For any
$y\in\partial_{\infty}(K\backslash G)$ and any $g\in G$ we have
$\lambda_{g\cdot x}(y)=\lambda_x(y\cdot g).$
\end{lemma}
\proof  Let $\gamma$ be a geodesic ray representing $y$. Using the
cocycle formula (\ref{eq:cocycle}) we compute
\begin{align*}
\lambda_{g\cdot x}(y) &= \lim_{t\to\infty}\phi_{g\cdot
x}(\gamma(t))= \lim_{t\to\infty}\frac{\Psi_{g\cdot
x}(\gamma(t))}{d(\gamma(t),x_0)}=
\lim_{t\to\infty}\frac{\Psi_{x}(\gamma(t)g)-\Psi_x(g)}{d(\gamma(t),x_0)}
\\
&=\lim_{t\to\infty}\frac{\Psi_{x}(\gamma(t)g)}{d(\gamma(t),x_0)}
=\lim_{t\to\infty}\frac{\Psi_{x}(\gamma(t)g)}{d(\gamma(t)g,x_0)}
=\lambda_x(y\cdot g),
\end{align*}
since $t\mapsto \gamma(t)g$ represents $y\cdot g$ and the quotient
$d(\gamma(t),x_0)/d(\gamma(t)g,x_0)$ converges to $1$ as
$t\to\infty$, because $d(\gamma(t),x_0)=d(\gamma(t)g,x_0g)$
converges to $\infty$ and by the triangular inequality
$|d(\gamma(t)g,x_0g)-d(\gamma(t)g,x_0)|\leq d(x_0g,x_0)$, which is
independent of $t$. \qed

\subsection{Some easy properties of maximal weights}

In the next two lemmata $x,x'$ denote points in $X$ and $s,s'$
denote elements in $\klie$. Recall that $\xi_s$ denotes the vector
field on $X$ generated by the infinitesimal action of $s$.

\begin{lemma}
\label{lemma:oposats-zero-fixe} If $\lambda(x;s)=\lambda(x;-s)=0$
then $\xi_s(x)=0$.
\end{lemma}
\proof By (\ref{eq:integral-t-pes}) $\lambda(x;s)=0$ implies that
$\la\mu(x),s\ra\leq 0$, and $\lambda(x;-s)=0$ implies
$\la\mu(x),-s\ra\leq 0$. Combining both inequalities we have
$\la\mu(x),s\ra=0$. Using again (\ref{eq:integral-t-pes}) and the
equality $\lambda(x;s)=\la\mu(x),s\ra$ we obtain $\xi_s(x)=0$.
\qed

\begin{lemma}
\label{lemma:max-weight-commuting} If $[s,s']=0$ then
$\lambda(x;s+s')=\lambda(x;s)+\lambda(x';s)$.
\end{lemma}
\proof We have
\begin{align*}
\lambda(x;s+s') &= \lim_{t\to\infty}\la \mu(e^{\imag t(s+s')}\cdot
x),s+s'\ra \\
&=\lim_{t\to\infty} \la \mu(e^{\imag ts}e^{\imag ts'}\cdot x),s\ra
+ \la \mu(e^{\imag ts}e^{\imag ts'}\cdot x),s'\ra \qquad\text{by
linearity and $[s,s']=0$} \\ &=\lim_{t\to\infty} \la \mu(e^{\imag
ts}\cdot x),s\ra + \la \mu(e^{\imag ts'}\cdot x),s'\ra
\qquad\text{by equivariance of $\mu$ and $[s,s']=0$} \\
&=\lambda(x;s)+\lambda(x;s').
\end{align*}
\qed

\section{Proof of Theorem \ref{thm:analytic-stab-zero-moment-map}}
\label{s:proof-theorem}

\subsection{Proofs of (1) and (2)}
\label{ss:proof-1-2} Statement (1) follows immediately from Lemma
\ref{lemma:lambda-equivariant}, and the fact that the action of
$G$ on $\partial_{\infty}(K\backslash G)$ extends the action by
isometries on $K\backslash G$ (so that the action on
$\partial_{\infty}(K\backslash G)$ of any element in $G$ sends
geodesically connected points to geodesically connected points).
(2) is well known, but we recall the argument for the sake of
completeness. If $G\cdot x\cap\mu^{-1}(0)$ contains two different
$K$-orbits, say $\oO_1,\oO_2\subset K$, then by Cartan's
decomposition we may find points $x_j\in\oO_j$ such that
$x_2=e^{\imag ts}\cdot x_1$. We then have
$\mu_s(x_2)=\mu_s(e^{\imag s}\cdot x)=\mu_s(x_1)=0$. By
(\ref{eq:gradient-mu-s}) we have
$$\mu_s(e^{\imag s}\cdot x)-\mu_s(x)=\int_0^1|\xi_s(e^{\imag\tau
s}\cdot x_1)|^2\,d\tau,$$ which implies that $\xi_s(e^{\imag\tau
s}\cdot x_1)$ vanishes for all $\tau\in[0,1]$, so that the action
of $\{e^{\imag ts}\mid t\in\RR\}$ fixes $x_1$. Consequently,
$x_2=x_1$, which implies $\oO_1=\oO_2$, a contradiction.

\subsection{Proof of (3)} Assume that $x$ is analytically
stable, so that for any $s\in S(\klie)$ we have $\lambda(x;s)>0$.
The usual argument to prove that $G\cdot
x\cap\mu^{-1}(0)\neq\emptyset$ is based on the identification
between the zeroes of the moment map and the critical values of
the integral of the moment map $\Psi_x$, and the fact that
analytic stability implies that $\Psi_x$ is proper. Instead, we
give here a topological argument. The condition $\lambda(x;s)>0$
implies that there is some $\tau_s$ such that if $t\geq \tau_s$
then $\la\mu(e^{\imag rs}\cdot x),s\ra>0$. Since the latter
function is continuous and $S(\klie)$ is compact, we may take some
$\tau$ working for any choice of $s$, namely, such that:
\begin{equation}
\label{eq:index-no-zero} \text{for any $t\geq\tau$ and any $s\in
S(\klie)$ we have $\la\mu(e^{\imag ts}\cdot x),s\ra>0$.}
\end{equation}
Denote by $\alpha:\klie^*\simeq\klie$ the isomorphism given by the
pairing $\la\cdot,\cdot\ra$. Property (\ref{eq:index-no-zero})
implies that the image of the map $f:S(\klie)\to\klie$ given
$f(s)=\alpha\circ\mu(e^{\imag \tau s}\cdot x)$ is contained in
$\klie\setminus\{0\}$, and furthermore there is a homotopy between
$f$ and the identity via maps from $S(\klie)$ to
$\klie\setminus\{0\}$. In other words, the index of $f$ around
$0\in\klie$ is nontrivial, and this implies that there is some $u$
inside the ball in $\klie$ with boundary $S(\klie)$ such that
$f(u)=0$, which is equivalent to $\mu(e^{\imag \tau u}\cdot x)=0$.
Now to prove that $G_x$ is finite is equivalent to proving that
$G_y$ is finite, where $y=e^{\imag\tau u}\cdot x$. Since
$\mu(y)=0$ and $\lambda(y;s)>0$ for any $s$, formula
(\ref{eq:integral-t-pes}) implies that for any $s\in\klie$ the
vector field $\xi_s$ is nonzero at $y$. Consequently the
stabilizer $K_y$ is finite. Finally, the condition $\mu(y)=0$
implies that $G_y$ is the complexification of $K_y$ (this is
proved by checking, using (\ref{eq:gradient-mu-s}), that if
$ke^{\imag u}$ fixes $y$, $k\in K$ and $u\in\klie$, then
$\xi_u(y)=0$ and $k\in K_x$, see \cite[Proposition 1.6]{S}).
Hence, $G_y$ is also finite.

The converse implication in (3) is almost immediate: if $y\in
G\cdot x\cap\mu^{-1}(0)$ and $G_x$ is finite, then $G_y$ is also
finite. This implies that $\xi_u(y)\neq 0$ for any $u\in
S(\klie)$, and now (\ref{eq:gradient-mu-s}) implies that
$\lambda(y;s)>0$, so $y$ is analytically stable. From (1) it now
follows that $x$ is also analytically stable.

\subsection{Proof of (4)} We first prove that
if $x$ is polystable then $G\cdot x\cap\mu^{-1}(0)\neq\emptyset$.
Since (3) has been proved, we only need to consider strictly
polystable points $x$ (namely, unstable polystable points). So let
$x\in X$ be such a point. Then one can choose $s\in S(\klie)$ such
that $\lambda_x(e_s)=0$ and such that $\dim T_{s'}\leq\dim T_s$
for any other $s'\in S(\klie)$ satisfying $\lambda_x(e_{s'})=0$.

Let $y=e_s$. Since $x$ is polystable, there exists some $y'\in
\partial_{\infty}(K\backslash G)$ which is geodesically connected to $y$
and such that $\lambda_x(y')=0$. Let $\gamma:\RR\to K\backslash G$
be a parameterized geodesic in $K\backslash G$ connecting $y$ and
$y'$, and assume that $\gamma(t)=e^{\imag ut}h$ for some
$u\in\klie$ and $h\in G$. By (1) the point $w=h\cdot x\in X$ is
polystable. If we set $u=s\cdot h^{-1}$ then we have $s'\cdot
h^{-1}=-u$, since the points $e_{s\cdot h^{-1}}$ and $e_{s'\cdot
h^{-1}}$ are connected by a geodesic passing through $x_0\in
K\backslash G$. By Lemma \ref{lemma:lambda-equivariant} we have
$\lambda_w(e_u)=\lambda_w(e_{-u})=0$. In other words,
$\lambda(w;u)=\lambda(w;-u)=0$. Then Lemma
\ref{lemma:oposats-zero-fixe} implies that $\xi_u(w)=0$. Hence the
group $\{\exp(ts)\mid t\in\RR\}\subset K$ fixes $w$, and by
continuity this implies that $T_u$ fixes $w$. So, if $\tlie_u$
denotes the Lie algebra of $T_u$ then for any $u'\in\tlie_u$ we
have $\xi_{u'}(w)=0$.

\begin{lemma}
\label{lemma:tot-tlie-u-es-zero} For any $u'\in\tlie_u$ we have
$\lambda(w;u')=0$.
\end{lemma}
\proof Since $w$ is polystable we have $\lambda(w;u')\geq 0$ and
$\lambda(w;-u')\geq 0$. Now, (\ref{eq:integral-t-pes}) together
with $\xi_{u'}(w)=0$ implies that
$\lambda(w;u')=\la\mu(w),u'\ra=-\la\mu(w),-u'\ra=-\lambda(w;-u')$.
\qed

Lemmata \ref{lemma:dim-torus-invariant} and
\ref{lemma:lambda-equivariant} imply that $u$ has the same
maximality property as $s$, namely
\begin{equation} \label{eq:u-maximal}
\text{$\dim T_{u'}\leq\dim T_u$ for any $u'\in S(\klie)$
satisfying $\lambda_w(e_{u'})=0$.}
\end{equation}

Let $K_u=\{k\in K\mid \Ad(k)(u)=u\}$ be the centralizer of $u$.
Then $T_u\subset K_u$ is obviously central. Let $K_0=K_u/T_u$ and
let $\klie_0$ be its Lie algebra. Consider the following maps:
\begin{enumerate}
\item the projection $\pi_u:\klie^*\to\klie_u^*$ induced by the
inclusion $\klie_u\subset\klie$, and \item the projection
$\pi_0:\klie_u^*\to\klie_0^*$ induced by any linear map
$\klie_0\to\klie_u$ which is a section of the projection
$\klie_u\to\klie_u/\tlie_u=\klie_0$ ($\pi_0$ is automatically a
morphism of Lie algebras because $\tlie_u$ is central in
$\klie_u$).
\end{enumerate}

Let $X_u\subset X$ be the set of points fixed by all elements of
$T_u$. Then $X_u$ is a Kaehler submanifold of $X$ and the group
$K_0$ acts on it by isometries. A moment map for this action,
$$\mu^{K_0}:X_u\to\klie_0^*,$$
can be obtained by composing
$\mu^{K_0}=\pi_0\circ\pi_u\circ\mu|_{X_u}$.

We claim that $w\in X_u$ is stable with respect to the action of
$K_0$. First of all we observe that for any $u'\in\tlie_u$ we have
$$\lambda^{K_0}(w;[u'])=\lambda^K(w;u'),$$ where on the
left hand side we consider the maximal weights of the action of
$K_0$ on $X_u$ and $[u']$ denotes the class in
$\klie_0=\klie_u/\tlie_u$ represented by $u'$, and on the right
hand side we consider the weights of the action of $K$ on $X$. It
follows that
$$\lambda^{K_0}(w;[u'])\geq 0$$
for any $u'\in \klie_u$. We claim that the latter inequality is
strict unless $[u']=0$. Indeed, if $[u']\neq 0$ and
$\lambda^{K_0}(w;[u'])=0$ then, letting $T\subset K$ be the torus
generated by $T_u$ and by the closure of $\{\exp(tu')\mid
t\in\RR\}$, we would have, by Lemma
\ref{lemma:max-weight-commuting} and arguing as in the proof of
Lemma \ref{lemma:tot-tlie-u-es-zero}, $\lambda(w;v)=0$ for any
$v\in\tlie=\Lie T$. Choosing $v$ in such a way that
$\{\exp(tv)\mid t\in\RR\}$ is dense in $T$ we would furthermore
have $\dim T_v>\dim T_u$, contradicting the maximality property
(\ref{eq:u-maximal}).

Hence $w$ is stable with respect to the action of $K_0$ on $X_u$,
so by (3) there exists some $h\in G_0$ such that $\mu^{K_0}(g\cdot
w)=0$. This immediately implies that $\mu^{K_u}(g\cdot w)=0$,
where $\mu^{K_u}=\pi_u\circ\mu|_{X_u}$ is the moment map for the
action of $K_u$ on $X_u$ (see Lemma
\ref{lemma:tot-tlie-u-es-zero}). We now prove that we also have
$\mu(g\cdot w)=0$. Let us denote for convenience $z=g\cdot w$.
Then $z$ is fixed by the action of $T_u$. Take a decomposition
$$\klie=\klie_u^*\oplus\bigoplus_{\alpha}\klie_{\alpha}$$
in irreducible representations of $T_u$, so that $\klie_u$ is the
trivial representation and each $\klie_{\alpha}$ is nontrivial.
This splitting induces a splitting of the dual vector space
$\klie^*$, and we let $\mu(z)=\mu_u(z)+\sum\mu_{\alpha}(z)$ be the
corresponding decomposition of $\mu$. We clearly have
$\mu_u(z)=\mu^{K_u}(z)$=0. Now, since $z$ is fixed by $T_u$, the
equivariance of the moment map implies that each $\mu_{\alpha}(z)$
is a $T_u$ invariant linear map $\tlie_{\alpha}\to\RR$. But each
$\tlie_{\alpha}$ is a nontrivial irreducible representation of
$T_u$, so the following lemma implies that $\mu_{\alpha}(z)=0$.

\begin{lemma}
Let $V$ be a finite dimensional vector space
and let $\Gamma\circlearrowright V$ be an irreducible nontrivial
linear action. Any $\Gamma$-invariant linear function $f:V\to\RR$
vanishes identically.
\end{lemma}
\proof Take any nonzero $v\in V$ which is not fixed by $\Gamma$.
Then the affine closure\footnote{If $X\subset V$, the affine
closure $\la X\ra_{\aff}\subset V$ is the set of finite sums
$\sum\lambda_ix_i$ with $\sum\lambda_i=1$ and $x_i\in X$.}
$\la\Gamma\cdot v\ra_{\aff}$ of $\Gamma\cdot v$ equal to $V$.
Indeed, $\la\Gamma\cdot v\ra_{\aff}$ is $\Gamma$-invariant and is
not a point, so if it were a proper subspace of $V$ then its
translate containing the origin would be a proper nonzero
invariant vector subspace of $V$, contradicting the irreducibility
of $\Gamma\circlearrowright V$. Since $f$ is $\Gamma$-invariant,
$f$ is constant on $\Gamma\cdot v$, and by linearity the
restriction of $f$ to $\la\Gamma\cdot v\ra_{\aff}$ is also
constant. Since $0\in \la\Gamma\cdot v\ra_{\aff}$ and $f$ is
linear, we must have $f=0$. \qed

Since $z\in G\cdot x$, we have proved that $G\cdot
x\cap\mu^{-1}(0)\neq\emptyset$.

The converse statement in (4) is almost immediate. Assume that
$G\cdot x\cap\mu^{-1}(0)\neq\emptyset$, and let $z\in G\cdot
x\cap\mu^{-1}(0)$. By statement (1) it suffices to prove that $z$
is polystable. Since $\mu(z)=0$, (\ref{eq:integral-t-pes}) implies
that $\lambda(z;s)\geq 0$ for any $s\in S(\klie)$, and also
$\lambda(z;s)=0$ if and only if $\xi_s(z)=0$. The latter implies
that $\lambda(z;s)=0$ if and only if $\lambda(z;-s)=0$. Since
$e_s,e_{-s}$ are always geodesically connected (see Example
\ref{example:s-minus-s}), it follows that $z$ is polystable.

It remains to prove that the stabilizer of polystable points is
reductive. Since we have proved that if $x$ is polystable then
$G\cdot x\cap\mu^{-1}(0)\neq\emptyset$, it suffices to prove that
if $\mu(z)=0$ then $G_z$ is reductive. This follows from the well
known observation that $G_z$ is the complexification of the
compact group $K_z=\{k\in K\mid k\cdot z=z\}$.

\section{Opposed elements in $\klie$ and geodesically connected
points in $\partial_{\infty}(K\backslash G)$}
\label{s:proof:lemma:opposed-geodesically-connected}

The main result of this section is the proof of Lemma
\ref{lemma:opposed-geodesically-connected}, which will be given in
the Section \ref{ss:lemma:opposed-geodesically-connected}. In
Sections \ref{ss:K-G-orbits-the-same} and \ref{ss:dense-orbit-OOO}
we state and prove some preliminary lemmata. Some of these results
are probably well known to experts, but we prove them in some
detail for the reader's convenience.

\subsection{The $K$-orbits and the $G$-orbits in
$\partial_{\infty}(K\backslash G)$ are the same}
\label{ss:K-G-orbits-the-same} Recall that for any $s\in S(\klie)$
the stabilizer of $e_s\in\partial_{\infty}(K\backslash G)$ is the
parabolic subgroup
\begin{equation}
\label{eq:def-parabolic} P_s=\{g\in G\mid e^{\imag ts}g e^{-\imag
ts}\text{ stays bounded as $t\to\infty$ }\}.
\end{equation}
As previously, we denote by $x_0 \in K\backslash G$ the class of
the identity element $1_G\in G$.

\begin{lemma}
\label{lemma:parabolic-transitive} For any $s\in S(\klie)$ the
action of $P_s$ on $K\backslash G$ (given by restricting the
action of $G$) is transitive.
\end{lemma}
\proof Take any $s\in S(\klie)$. Since $K\backslash G$ is
connected, to prove the lemma it suffices to check that $x_0\cdot
P_s$ is open and closed in $K\backslash G$. With the aim of
proving that $x_0\cdot P_s$ is open, let us check that
$\klie+\plie_s=\glie$, where $\plie_s\subset\glie$ is the Lie
algebra of $P_s$. The endomorphism $\ad(s)\in\End\glie$ is
semisimple because it preserves the extension to $\glie$ of the
biinvariant scalar product in $\klie$. Hence we may consider the
decomposition in eigenspaces $\glie=\bigoplus\glie_{\lambda}$ of
the action of $\imag\ad(s)$ on $\glie$, where each $\lambda$ is
real and $\imag\ad(s)$ acts on $\glie_{\lambda}$ as multiplication
by $\lambda$. It follows from (\ref{eq:def-parabolic}) that
$\plie_s=\bigoplus_{\lambda\leq 0}\glie_{\lambda}$. Let
$c:\glie\to\glie$ denote the conjugation map given by the
identification $\glie\simeq\klie\otimes_{\RR}\CC$, so that
$\klie\subset\glie$ is the fixed point set of $c$. Since
$c([a,b])=[c(a),c(b)]$ for any $a,b\in\glie$, we have $c\circ
\imag\ad(s)\circ c=-\imag\ad(s)$, which implies that $c$ induces
isomorphisms $\glie_{\lambda}\simeq\glie_{-\lambda}$ for each
$\lambda$. Since $\klie$ is the fixed point set of $c$, for any
nonzero $\lambda$ the intersection
$\klie\cap(\glie_\lambda\oplus\glie_{-\lambda})$ is equal to the
graph of $c:\glie_{-\lambda}\to\glie_{\lambda}$. Combining this
with the fact that $\plie_s=\bigoplus_{\lambda\leq
0}\glie_{\lambda}$, we deduce that $\klie+\plie_s=\glie$. This
implies, by the inverse function theorem, that any $g\in G$
sufficiently near $1_G$ can be written as $g=k\cdot p$ for some
$k\in K$ and $p\in P_s$, which means that $x_0\cdot P_s$ contains
a neighborhood of $x_0$. Since $P_s$ acts on $K\backslash G$ by
homeomorphisms, this implies that $x_0\cdot P_s$ contains a
neighborhood of any of its points, so it is open.

Let $\alpha:P_s\to K\backslash G$ be the map $p\mapsto
\alpha(p)=x_0\cdot p$. Then $\alpha$ is the restriction of the
quotient map $G\to K\backslash G$, which is proper because $K$ is
compact. Since $P_s\subset G$ is closed, it follows that $\alpha$
is also proper, so the intersection of $\alpha(P_s)$ with any
compact subset of $K\backslash G$ is closed. Since $K\backslash G$
is locally compact, it follows that $x_0\cdot P_s=\alpha(P_s)$ is
closed.\qed

\begin{lemma}
\label{lemma:orbits-in-boundary} If
$y,y'\in\partial_{\infty}(K\backslash G)$ satisfy $y'=y\cdot g$
for some $g\in G$, then there exists some $k\in K$ such that
$y'=y\cdot k$.
\end{lemma}
\proof Assume that $y'=y\cdot g=e_s$ and $z=x_0\cdot g$. By Lemma
\ref{lemma:parabolic-transitive} there exists some $p\in P_s$ such
that $z\cdot p=x_0$. Hence $k:=gp$ satisfies $y'=y\cdot k$ and
$x_0\cdot k=x_0$, which implies that $k\in K$. \qed

The previous lemma implies that for any $s\in\klie$ of unit norm
there is a right action of $G$ on the adjoint orbit
$\oO_s=\Ad(K)\cdot s\subset\klie$. Indeed, via the map $S(\klie)
\ni u\mapsto e_u\in \partial_{\infty}(K\backslash G)$ the action
of $K$ on the boundary $\partial_{\infty}(K\backslash G)$
corresponds to the adjoint action on $S(\klie)$. Since the
$K$-orbits in $\partial_{\infty}(K\backslash G)$ are equal to the
$G$-orbits, for any $s\in S(\klie)$ we can identify $\oO_s$ with
one of the $G$-orbits. And since the stabilizer of $s$ is $P_s$,
we obtain a natural identification $\oO_s\simeq P_s\backslash G$.

\subsection{The dense orbit of the action of $P_{-s}$ on $\oO_s$}
\label{ss:dense-orbit-OOO}

Let $\oO_{s}^*\subset\oO_{s}$ denote the set of elements which are
opposed to $-s$.

\begin{lemma}
\label{lemma:O-star-connected} The set $\oO_s^*\subset\oO_s$ is
open, dense, and connected.
\end{lemma}
\proof To prove the lemma we check that $\oO_s$ carries a
structure of complex connected manifold with respect to which
$\oO_s\setminus\oO_s^*\subset\oO_s$ is an analytic subvariety of
dimension $<\dim\oO_s$. Since $K$ is connected, $\oO_s=\Ad(K)(s)$
is also connected. Let $\oO_{\ad(s)}\subset\End\glie$ be the
adjoint orbit of $\ad(s)$ under the action of the vector space
automorphisms of $\glie$. Let $\glie=\bigoplus\glie_{\lambda}$ be
the eigenspace decomposition of the action of $\imag\ad(s)$ on
$\glie$, as in the proof of Lemma
\ref{lemma:parabolic-transitive}. Let $\fF$ be the set of growing
filtrations $(W^{\mu})_{\mu\in\RR}$ of complex subspaces of
$\glie$ satisfying $\dim
W^{\mu}=\sum_{\lambda\leq\mu}\dim\glie_{\lambda}$. The set $\fF$
(which is an example of flag variety) carries a natural structure
of complex manifold, and the map $w:\oO_{\ad(s)}\to\fF$ which
sends $u\in\oO_{\ad(s)}$ to the filtration $(W_u^{\mu})$ with
$W_u^{\mu}=\bigoplus_{\lambda\leq\mu}\Ker(\imag u-\lambda\Id)$ is
clearly a diffeomorphism. Let $f:\oO_{s}\to \oO_{\ad(s)}$ be the
restriction of $\ad:\klie\to\End\glie$. We claim that the map
$$\phi=w\circ f:\oO_s\to\fF$$ is an immersion and that for any
$v\in\oO_s$ the image $d\phi(T_v\oO_s)\subset T_{\phi(v)}\fF$ is a
complex subspace, so that there is a unique structure of complex
manifold on $\oO_s$ with respect to which $\phi$ is holomorphic.
Since for any $h\in K$ we have
\begin{equation}
\label{eq:phi-equivariant}
\phi(\Ad(h)(s))=w(\Ad(h)\ad(s)\Ad(h)^{-1})=\Ad(h)\phi(s)
\end{equation}
and the map $\Ad(h):\fF\to\fF$ is a biholomorphism, to prove the
claim it suffices to check that $d\phi(s):T_s\oO_s\to
T_{\phi(s)}\fF$ is an injection and that its image is invariant
under multiplication by $\imag$.

Proving that $d\phi(s)$ is injective is equivalent to proving that
$df(s):T_s\oO_s\to T_{\ad(s)}\oO_{\ad(s)}$ is injective, because
$w$ is a diffeomorphism. We have $T_s\oO_s=\{[a,s]\mid
a\in\klie\}$. Assume that $df(s)([a,s])=0$. Since $f$ is the
restriction of the map $\ad:\klie\to\End\glie$, which is linear,
we deduce from the assumption that $\ad([a,s])=0$. Write $a=\sum
a_{\lambda}$, where $a_{\lambda}\in\glie_{\lambda}$. Then
$[a,s]=-\sum\lambda a_{\lambda}$, and similarly
$0=\ad([a,s])(s)=[[a,s],s]=\sum\lambda^2a_{\lambda},$ so that
$a_{\lambda}=0$ for each $\lambda\neq 0$. But then we have
$[a,s]=-\sum\lambda a_{\lambda}=0$. Hence we have proved that
$df(s)$ is injective.

We now check that $d\phi(s)(T_s\oO_s)\subset T_{\phi(s)}\fF$ is a
complex subspace. Here we give a direct argument but an
alternative and more intrinsic proof of this result may be given
using Lemma \ref{lemma:phi-G-equivariant} below. Remark that
$$T_{\phi(s)}\fF=\bigoplus_{\lambda<\mu}
\Hom(\glie_{\lambda},\glie_{\mu}).$$ Take any $a\in\klie$. By
(\ref{eq:phi-equivariant}), $\phi(s)$ sends $[a,s]$ to the
projection of
$\ad(a)\in\bigoplus_{\lambda,\mu}\Hom(\glie_{\lambda},\glie_{\mu})$
to $T_{\phi(s)}$. So if we decompose $a=\sum a_{\lambda}$ as
before, then the piece of $d\phi(s)(a)$ in
$\Hom(\glie_{\lambda},\glie_{\mu})$ is $\ad(a_{\mu-\lambda})$ (in
particular $d\phi(s)(a)$ only depends on
$\sum_{\lambda>0}a_{\lambda}$). Denote by $c:\glie\to\glie$ the
conjugation coming from identifying $\glie=\klie\otimes_{\RR}\CC$,
as in the proof of Lemma \ref{lemma:parabolic-transitive}. Since
$\klie$ is the fixed point set of $c$ and
$c(\glie_{-\lambda})=\glie_{\lambda}$, it follows from $a\in\klie$
that
$$a'=\sum_{\lambda<0}-\imag a_{\lambda}+a_0+\sum_{\lambda>0}\imag
a_{\lambda}\quad\text{also belongs to }\klie.$$ But by the
previous observations we have $d\phi(s)([a',s])=\imag
d\phi(s)([a,s])$. Hence the image of $d\phi(s)$ is invariant under
multiplication by $\imag$, which is what we wanted to prove.

To finish the proof of the lemma, let $\fF^*$ be the set of
filtrations $(W^{\mu})_{\mu}\in\fF$ such that
$$\glie=\bigoplus_{\lambda} \left( W^{\lambda}\cap
(\bigoplus_{\nu\geq\lambda}\glie_{\nu})\right),$$ where the sum
runs over the spectrum of $\imag\ad(s)$. It is straightforward to
check that $\fF\setminus\fF^*$ is an analytic subvariety of $\fF$.
Since $\phi:\oO_s\to\fF$ is a holomorphic map and
$\oO_s^*=\phi^{-1}(\fF^*)$, we deduce that $\oO_s\setminus\oO_s^*$
is an analytic subvariety of $\oO_s$. Finally, since $\oO_s^*$ is
nonempty (it contains $s$, for example), $\oO_s\setminus\oO_s^*$
is not equal to $\oO_s$, and since $\oO_s$ is connected this
implies that $\dim \oO_s\setminus\oO_s^*<\dim\oO_s$. \qed

\begin{lemma}
\label{lemma:phi-G-equivariant} Let $\phi:\oO_s\to\fF$ be the map
defined in the proof of Lemma \ref{lemma:O-star-connected}. For
any $u\in\oO_s$ and $g\in G$ we have $\phi(u\cdot
g)=\Ad(g^{-1})\phi(u)$.
\end{lemma}
\proof We first state a general result relating morphisms between
groups and morphisms between boundaries of the corresponding
symmetric spaces. Let $\rho:K\to K'$ be a morphism of compact
connected Lie groups, and denote by the same symbol $\rho:G\to G'$
the induced map between the complexifications. Let $x_0=[1_G]\in
K\backslash G$ and $x_0'=[1_{G'}]\in K'\backslash G'$ be the
classes of the identity elements. There is a unique map
$r:K\backslash G\to K'\backslash G'$ satisfying $r(x_0\cdot
g)=x_0'\cdot \rho(g)$ for any $g\in G$. Choosing biinvariant
metrics on the Lie algebras of $K$ and $K'$ and taking the induced
Riemannian structures on $K\backslash G$ and $K'\backslash G'$,
the map $r$ is Lipschitz. Furthermore, $r$ sends geodesic rays in
$K\backslash G$ either to constant maps or to geodesic rays in
$K'\backslash G'$. More precisely, if $t\mapsto
\gamma(t)=[e^{\imag ts}g]\in K\backslash G$ is a geodesic ray,
then: $t\mapsto r\circ\gamma(t)$ is a geodesic ray in
$K'\backslash G'$ unless $d\rho(s)=0$, in which case we obtain a
constant map. Since $r$ is Lipschitz, given any pair of equivalent
geodesic rays $\gamma_0\sim\gamma_1$ either both $r\circ\gamma_0$
and $r\circ\gamma_1$ are geodesic rays or both are constant maps.
So one may define the set $\partial_{\infty}(K\setminus G)^*$ of
boundary points corresponding to geodesic rays which are mapped to
geodesic rays, and then $r$ induces a continuous map
$$r:\partial_{\infty}(K\backslash G)^*\to \partial_{\infty}(K'\backslash
G').$$ It is clear that $\partial_{\infty}(K\backslash G)^*\subset
\partial_{\infty}(K\backslash G)$ is $G$-invariant and that $r$ is
equivariant, in the sense that for any $y\in
\partial_{\infty}(K\backslash G)^*$ and $g\in G$ we have $r(y\cdot
g)=r(y)\cdot\rho(g)$.

Consider the maps $f:\oO_s\to\oO_{\ad(s)}$ and
$w:\oO_{\ad(s)}\to\fF$ given in the proof of Lemma
\ref{lemma:O-star-connected}. We apply the previous observations
to the case in which $K'$ is the set of vector space automorphisms
of $\glie$ preserving the Hermitian product induced by the
biinvariant metric on $\klie$. The complexification of $K'$ is the
group $G'$ of all automorphisms of $\glie$, and we may take as a
morphism $\rho:K\to K'$ the adjoint representation:
$\rho(k)=\Ad(k)$. The conclusion is that $f(u\cdot
g)=f(u)\cdot\Ad(g)$ for any $u\in \oO_s$ and $g\in G$. Finally, by
the results in Section \ref{ss:boundary-unitary} we also have
$w(u\cdot g)=\Ad(g^{-1})w(g)$ for any $g\in G'$ and $u\in\klie'$
of unit norm. This finishes the proof of the lemma. \qed

It follows from the previous lemma that the set
$\oO_s^*\subset\oO_s$ is $G$-invariant. We next prove that
$\oO_s^*$ is an orbit of the induced action of $P_{-s}\subset G$
on $\oO_s$.

\begin{lemma}
\label{lemma:parabolic-transitive-on-opposed} The action of
$P_{-s}$ on $\oO_{s}^*$ is transitive.
\end{lemma}
\proof Let $\ulie=\bigoplus_{\lambda>0}\glie_{\lambda}$. This is
the Lie algebra of the biggest unipotent subgroup of $P_{-s}$.
Consider the map $e:\ulie\to \oO_s^*$ defined as $e(u)=s\cdot
e^u$. We are going to prove that the image of $e$ is $\oO_s^*$.
Since by Lemma \ref{lemma:O-star-connected} $\oO_s^*$ is
connected, it suffices to prove that $e(\ulie)$ is open and closed
in $\oO_s^*$. From $\plie_s=\bigoplus_{\lambda\leq
0}\glie_{\lambda}$ we deduce that $\glie=\plie_s\oplus\ulie$ which
implies, by the implicit function theorem, that any $g\in G$
sufficiently close to $1_G$ can be written as $g=pe^{u}$, where
$p\in P_s$ and $u\in \ulie$. Hence $e(\ulie)$ is open in
$P_s\backslash G$ (use the same arguments as in the proof of Lemma
\ref{lemma:parabolic-transitive}). Let the map $\phi:\oO_s\to\fF$
and the subset $\fF^*\subset\fF$ be those defined in the proof of
Lemma \ref{lemma:O-star-connected}. We then have
$\phi(\oO_s^*)\subset\fF^*$. There is a biholomorphism
$$\gamma:\fF^*\to\bigoplus_{\lambda<\mu}
\Hom(\glie_{\lambda},\glie_{\mu}),$$ characterized by the property
that $\gamma^{-1}$ sends $\delta=(\delta_{\lambda\nu})\in
\bigoplus_{\lambda<\mu} \Hom(\glie_{\lambda},\glie_{\mu})$ to the
filtration $(W^{\mu}(\delta))_{\mu}$ in which
$W^{\mu}(\delta)=\bigoplus_{\xi\leq\mu}\Graph(\delta_{\xi})$,
where $\delta_{\xi}= \sum_{\mu\geq\xi}\delta_{\xi\mu}:\glie_{\xi}
\to\bigoplus_{\mu\geq\xi} \glie_{\mu}$. Now, to check that
$e(\ulie)\subset\oO_s^*$ is closed it suffices to prove, similarly
to Lemma \ref{lemma:parabolic-transitive}, that the map $f$
defined as the following composition is proper:
$$f:\ulie\stackrel{e}{\longrightarrow}\oO_s^*
\stackrel{\phi}{\longrightarrow}\fF^*
\stackrel{\gamma}{\longrightarrow} \bigoplus_{\lambda<\mu}
\Hom(\glie_{\lambda},\glie_{\mu}).$$

Let $0<\lambda_1<\dots<\lambda_r$ be the positive eigenvalues of
$\imag\ad(s)$. For any $u=\sum_{\lambda_i}\in\ulie$ write
$u_i:=u_{\lambda_i}$, and let $f(u)_j$ denote the component of
$f(u)$ in $\Hom(\glie_0,\glie_{\lambda_j})$. We deduce from the
definitions that $f(u)=\exp(-\ad(u))-1$ for any $u\in\ulie$. Since
by the Jacobi identity
$[\glie_{\lambda},\glie_{\mu}]\subset\glie_{\lambda+\mu}$, and the
decomposition $\glie=\bigoplus\glie_{\lambda}$ is finite, there
exist polynomials $P_j$ such that $P_j(0,\dots,0)=0$ and
\begin{equation}
\label{eq:f-escalat}
f(u)_j=\ad(u_j)_j+P_j(\ad(u_1),\dots,\ad(u_{j-1})),
\end{equation}
where $\ad(u_j)_j$ denotes the piece of $\ad(u_j)\in\End\glie$ in
$\Hom(\glie_0,\glie_{\lambda_j})$.

Consider the Hermitian norm on $\glie$ induced by the biinvariant
norm on $\klie$ and define, for any $\alpha\in\End\glie$,
$|\alpha|=\sup |\alpha(v)|/|v|$, where the supremum runs over the
set of nonzero $v\in\glie$. Given any $u\in\ulie$, we have for
each $j$:
\begin{equation}
\label{eq:norm-ad-u-j} |\ad(u_j)_j|\geq
|\ad(u_j)(s)|/|s|=\lambda_j|u_j|.
\end{equation}
Let $f(u)_j$ be the piece of $f(u)$ in
$\Hom(\glie_0,\glie_{\lambda_j})$. Then we have
$f(u)_{\lambda_j}=\ad(u_{\lambda_j})$, so
$$|f(u)|\geq |f(u)_{\lambda_j}|\geq |f(u)_{\lambda_j}(s)|/|s|=
|[u_{\lambda_j},s]|/|s|=\lambda_j|u_{\lambda_j}|.$$ There exist
polynomials $p_j\in\CC[t]$ vanishing at $t=0$ such that
$$|P_j(\ad(u_1),\dots,\ad(u_{j-1}))|\leq
p_j(|u_1|+\dots+|u_{j-1}|)$$ for each $j$ and $u_1,\dots,u_{j-1}$.
Since $p_1(0)=\dots=p_r(0)=0$, there exists an $\epsilon>0$ such
that the following system of inequalities
\begin{align*}
t_1 &< \epsilon(t_1+\dots+t_r) \\
t_2 &< \epsilon(t_1+\dots+t_r) + 2\lambda_2^{-1}p_2(t_1) \\
& \vdots \\
t_r &< \epsilon(t_1+\dots+t_r)+
2\lambda_r^{-1}p_r(t_1+\dots+t_{r-1})
\end{align*}
has no solution $(t_1,\dots,t_r)$ satisfying $t_j\geq 0$ for each
$j$.  Let us prove that for any $u\in\ulie$ we have $|f(u)|\geq
\lambda_1\epsilon |u|/2$, which clearly implies that $f$ is
proper. Define $t_j=|u_j|$ for each $j$. By the choice of
$\epsilon$ at least one of the previous inequalities does not
hold, say the $j$-th one. Then we have (setting $p_1=0$ and
$P_1=0$ when $j=1$)
$$|u_j|\geq \epsilon(|u_1|+\dots+|u_r|)+
2\lambda_j^{-1}p_j(|u_1|+\dots+|u_{j-1}|)$$ which implies, using
(\ref{eq:norm-ad-u-j}) and the definition of $p_j$,
$$|\ad(u_j)_j|\geq\lambda_j|u_j|\geq 2p_j(|u_1|+\dots+|u_{j-1}|)
\geq 2|P_j(\ad(u_1),\dots,\ad(u_{j-1}))|.$$ Combining this with
(\ref{eq:f-escalat}) we obtain
\begin{align*}
|f(u)| & \geq |f(u)_j| \geq |\ad(u_j)_j|/2 \geq \lambda_j|u_j|/2
\\
&\geq \lambda_j\epsilon(|u_1|+\dots+|u_r|)/2 \\
&\geq \lambda_j\epsilon|u|/2\geq\lambda_1\epsilon|u|/2.
\end{align*}
This finishes the proof of the lemma. \qed

\subsection{Proof of Lemma
\ref{lemma:opposed-geodesically-connected}}
\label{ss:lemma:opposed-geodesically-connected} Assume that
$e_u,e_v$ are geodesically connected, and let $\gamma:\RR\to
K\backslash G$ be a geodesic such that $\gamma(t)\to e_u$ when
$t\to\infty$ and $\gamma(t)\to e_v$ when $t\to-\infty$. By Lemma
\ref{lemma:parabolic-transitive} there exists some $h\in P_u$ such
that $\gamma(0)\cdot h=x_0$. Since $\gamma\cdot h$ is a geodesic
passing through $x_0$ at time $0$, it is of the form $\gamma\cdot
h(t)=e^{\imag ts}$, and since $\gamma\cdot h(t)$ converges to
$e_u$ as $t\to \infty$, we have $s\cdot h=s=u$. Then $v\cdot
h=-u$. By Lemma \ref{lemma:orbits-in-boundary}, $v$ belongs to the
adjoint orbit $\oO_{-u}\subset\klie$. Obviously the endomorphisms
$\ad(u),\ad(-u)\in\End\glie$ are opposed (in the sense specified
in the Introduction), and Lemma \ref{lemma:phi-G-equivariant}
implies that $\ad(u),\ad(v)$ are opposed as well. Hence, $u$ and
$v$ are opposed.

Conversely, assume that $u$ and $v$ are opposed. Then
$v\in\oO_{-u}^*$. By Lemma
\ref{lemma:parabolic-transitive-on-opposed} the action of $P_u$ on
$\oO_{-u}^*$ is transitive, so there exists some $h\in P_u$ such
that $v\cdot h=-u$. Then the geodesic $\gamma(t)=[e^{\imag
tu}h^{-1}]$ satisfies $\gamma(t)\to e_u$ when $t\to\infty$ and
$\gamma(t)\to e_v$ when $t\to\infty$, so $e_u$ and $e_v$ are
geodesically connected. This finishes the proof of the Lemma.

\end{document}